\documentclass[12pt]{amsart}

\usepackage{amsmath,amsthm,amscd,euscript,amssymb,graphicx}
\def \r{\mathbb R}

\DeclareMathOperator{\cri}{Cr}
\DeclareMathOperator{\ind}{ind}

\newtheorem{theorem}{Theorem}[section]

\theoremstyle{remark}
\newtheorem{remark}[theorem]{Remark}
\theoremstyle{definition}
\newtheorem{definition}[theorem]{Definition}
\newtheorem{example}[theorem]{Example}

\title[Tensegrities on the space of generic functions]{Tensegrities on the space of generic functions}
\author{Oleg Karpenkov}
\date{25 May 2019}

\begin{document}
\maketitle
\input{epsf}

\section*{Introduction}

In this small note we introduce a notion of self-stresses on the set 
functions in two variables with generic critical points.
The notion naturally comes from a rather exotic representation of classical Maxwell frameworks 
in terms of differential forms.

For the sake of clarity we work in the two-dimensional case only. 
However all the definitions for the higher dimensional case are straightforward.

\section{Preliminaries}

\subsection{Classical definition of tensegrity}

For completeness of the story we start with the classical approach introduced 
in~\cite{Max} by J.~C.~Maxwell in 1864. We refer any interested in rigidity 
and flexibility questions to~\cite{Con,SJS2018}.

We use the following slightly modified definition of tensegrity from~\cite{Guz1}.

\begin{definition}
Let $G=(V,E)$ be an arbitrary graph on $n$ vertices.
\itemize
\item{A {\it framework} $G(P)$ in the plane is a map $f=(f_v,f_e)$:
$$
f_v: V \to \r^2, \qquad f_e: E \to S^{1},
$$
such that for every edge $v_iv_j$ the vector
$f_v(v_i)f_v(v_j)$ is a multiple of $f_e(v_iv_j)$.
}

\item{A {\it stress} $w$ on a framework is an assignment of real scalars
$w_{i,j}=w_{j,i}$ (called {\it tensions}) to its edges.}

\item{A stress $w$ is called a {\it self-stress} if, in addition, the
following equilibrium condition is fulfilled at every vertex
$p_i$:
$$
\sum\limits_{\{j|j\ne i\}} w_{i,j}e_{ij}=0.
$$
}
\item{A pair $(G(P),w)$ is a {\it tensegrity} if $w$ is a self-stress  for the framework $G(P)$.}
\end{definition}


\subsection{Tensegrities and exterior forms}

In this section we recall a rather exotic interpretation of two-dimensional tensegrities 
as a collection of 2-forms in $\r^3$ with certain relations.
It is rather in common with projective approach discussed by I.~Izmestiev in~\cite{Izm2017}.

\vspace{2mm}

Consider  a tensegrity $(G(P),w)$ with $P=(P_1,\ldots, P_n)$.
For every point $P_i=(x_i,y_i)$ we associate a 1-form in $\r^3$:
$$
dP_i:=x_idx+y_idy+dz.
$$
(one can say that $dz$ is a normalization factor to extract tensions.)

For every edge $E_iE_j$ we consider a 2-form:
$$
dP_i\wedge dP_j.
$$
It turns out that self-stressability conditions is precisely equivalent to 
$$
\sum\limits_{\{j|j\ne i\}} w_{i,j}dP_j\wedge dP_i=0.
$$
So any framework in tensegrity can be defined simply by a collection of 
decomposable 2-forms in $\r^3$ (which we denote as $G(dP)$) and a self-stress $w$ as before, 
We denote it by $(G(dP),w)$.

\begin{remark}
This definition perfectly suits the ``meet'' and ``join'' Cayley algebra expressions arising with
the description of existence conditions of tensegrities 
(see, e.g., in~\cite{WW,DKS2010,Kar2018}).
It also rather straightforwardly provides projective invariance of tensegrity existence.
\end{remark}


\section{Case of generic functions}

One of the mysterious questions related the notion of $(G(dP),w)$ is as follows: {\it what is a natural generalizations
of the tensegrity to the case of decomposable differentiable 2-forms $($not-necessarily with constant coefficients$)$?}

\vspace{1mm}

The aim of this section is to give a partial answer to this question for differential forms whose factors are of type
$$
df+dz, 
$$
where $f=f(x,y)$ is a function of two variables.

\vspace{2mm}

Let us first give a definition of tensegrity.
Secondly we show a geometric interpretation and link it to the classical case.

\subsection{Main definitions}
Let $F=(f_1,\ldots, f_n)$. Denote by $dF$ the collection of forms
$$
dF_1=df_1+dz, \quad dF_2=df_2+dz, \quad\ldots,\quad dF_n=df_n+dz.
$$

\begin{definition}
Let $F$ be a collection of functions $F=(f_1,\ldots, f_n)$ with finitely many critical points.
A {\it tensegrity} $(G(dF),w)$ is a triple: a graph $(G,F,w)$,
where functions $f_i$ are associated with vertices of a graph and edges are associated with stresses $w_{i,j}$.
\end{definition}

For a function $f$ denote the set of its critical points by $\cri(f)$;
the index of a critical point $P$ is denoted by $\ind(P)$.

\begin{definition}
A self-stress condition on $(G(dF),w)$ at a function $F_i$
$$
\sum\limits_{P_{i,k}\in \cri(f_i)} (-1)^{\ind(P_{i,k})}
\bigg(
\sum\limits_{\{j|j\ne i\}} w_{i,j}dF_j(P_{i,k})\wedge dF_i(P_{i,k})
\bigg)=0.
$$
\end{definition}

\begin{remark}
It might be also useful to consider critical points separately (say if this increases durability for certain overconstrained system).
\end{remark}

\begin{remark}
Recall that at a critical point of $F_i$
$$
dF_i\wedge dF_j=dz\wedge df_j.
$$
So one can replace every 2-form $dF_i\wedge dF_j$ in the equilibrium condition simply by $df_j$.
\end{remark}

\subsection{Geometric discussions}

\begin{figure}
\epsfbox{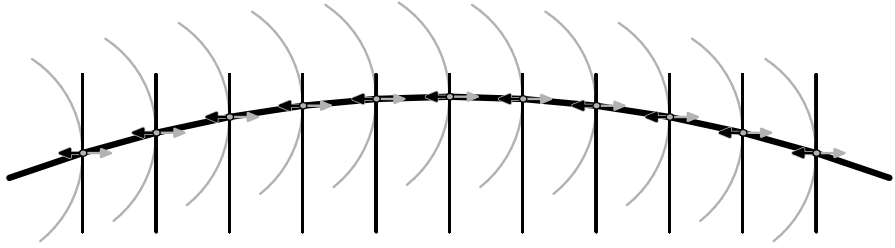}		
		\caption{Level sets of two functions (gray and black), and the line where their gradients have the same direction.}
		\label{f.1}
\end{figure}

{\it Lines of forces} are precisely the points when the total force $dF_i\wedge dF_j=0$
(see Figure~\ref{f.1}).
Lines of forces are defined by the equation 
$$
df_i\wedge df_j=0.
$$

It is clear that lines of forces are connecting critical points of $F_i$ (including critical points at infinity) to critical points $F_j$
usually by a graph rather than by a line. 
A possible picture for a splitting of a force line is as follows:
$$
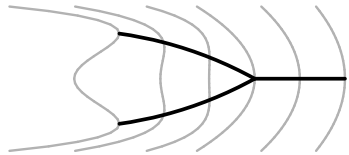
\epsfbox{f.2}
$$

\begin{remark}
One might consider the classical theory of tensegrities as follows:
For each point $P_i=(x_i,y_i)$ consider a function  
functions 
$$
f_i=(x-x_i)^2+(y-y_i)^2.
$$
Then (possible after a simple rescaling of stresses) one has a classical tensegrity.

This works also for case of point-hyperplane frameworks introduced recently in~\cite{EJN}, 
where hyperplanes are defined as linear functions.
\end{remark}

In some sense the proposed techniques can be considered as a deformation of a classical tensegrity.

\begin{figure}
\epsfbox{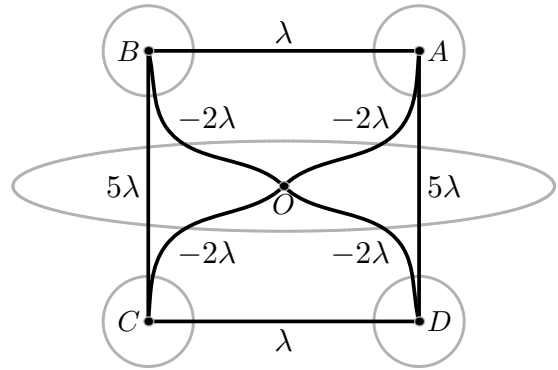}
		\caption{Level sets, compact lines of forces, and stresses.}
		\label{f.3}
\end{figure}

\begin{example}
First we start with a graph $G$ on 5 vertices:
$$
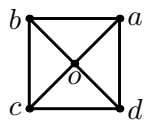
\epsfbox{f.4}
$$
Let us consider the following 5 functions corresponding to vertices of graphs:
$$
\begin{array}{l}
a:\quad f_1=(x-3)^2+(y-3)^2;\\
b:\quad f_2=(x+3)^2+(y-3)^2;\\
c:\quad f_3=(x+3)^2+(y+3)^2;\\
d:\quad f_4=(x-3)^2+(y+3)^2;\\
o:\quad f_5=x^2+5y^2.
\end{array}
$$ 
Then the line of forces and the corresponding stresses are as on Figure~\ref{f.3}.
They are defined up to a choice of a real parameter $\lambda$.
Here grey curves are level sets; black curves are compact lines of forces between critical points;
the numbers indicate the stresses on edges.
The critical points of functions are marked by the corresponding capital letters.    
\end{example}

\begin{remark}
Finally we would like to admit that 
the situation in three and higher dimensional cases repeats the two-dimensional case discussed above.
\end{remark}

\begin{remark}
In three dimensional case consider the following functions (potentials):
$$
f_{a,b,c}=\frac{k_e}{(x-a)^2+(y-b)^2+(z-c)^2}
$$
(here $k_e$ is the Coulomb constant) and take the unit stresses. 
Then we arrive to classical Coulomb situation for points with unit charges in three-space.

\end{remark}


\end{document}